\documentclass[12pt]{amsart}
\usepackage{amsfonts,amssymb,amsmath,amscd,amsthm}
\usepackage{graphicx} 
\usepackage{subfigure}
\usepackage[all]{xy}

\usepackage[colorlinks,allcolors=blue]{hyperref} 

\newtheorem{thm}{Theorem}
\newtheorem{lem}[thm]{Lemma}
\newtheorem{prop}[thm]{Proposition}

\newtheorem{cor}[thm]{Corollary}

\theoremstyle{definition}
\newtheorem{dfn}[thm]{Definition}
\newtheorem{ex}[thm]{Example}
\newtheorem{rmk}[thm]{Remark}

\newtheorem{question}[thm]{Question}

\numberwithin{thm}{section}
\numberwithin{equation}{section}

\newcommand{\Proof}{\noindent {\it Proof}.\ \ }

\newcommand{\Ext}{\operatorname{Ext}}
\newcommand{\Spec}{\operatorname{Spec}}

\newcommand{\ob}{\operatorname{ob}}
\newcommand{\Hilb}{\operatorname{Hilb}}
\newcommand{\HF}{\operatorname{HF}}
\newcommand{\red}{\operatorname{red}}

\renewcommand{\div}{\operatorname{div}}

\newcommand{\Pic}{\operatorname{Pic}}

\newcommand{\Bl}{\operatorname{Bl}}

\newcommand{\car}{\operatorname{char}}

\newcommand{\Fix}{\operatorname{Fix}}
\newcommand{\rank}{\operatorname{rank}}
\newcommand{\Sym}{\operatorname{Sym}}

\newcommand{\codim}[1]{\operatorname{codim}}

\renewcommand{\labelenumi}{{\rm (\arabic{enumi})}}

\title{Obstructions to deforming curves on an Enriques-Fano $3$-fold}
\author{Hirokazu Nasu}
\address{
Department of Mathematical Sciences,
Tokai University,
4-1-1 Kitakaname, Hiratsuka, 
Kanagawa 259-1292, JAPAN}
\email{nasu@tokai-u.jp}
\dedicatory{Dedicated to Professor~Shigeru~Mukai on the~occasion of his~65-th~birthday}

\begin{document}

\begin{abstract}

We study the deformations of a curve $C$ on 
an Enriques-Fano $3$-fold $X \subset \mathbb P^n$, assuming that $C$
is contained in a smooth hyperplane section $S \subset X$,
that is a smooth Enriques surface in $X$.
We give a sufficient condition for $C$ to be (un)obstructed in $X$,
in terms of half pencils and $(-2)$-curves on $S$.
Let $\Hilb^{sc} X$ denote the Hilbert scheme of smooth connected curves in $X$.
By using the Hilbert-flag scheme of $X$,
we also compute the dimension of $\Hilb^{sc} X$ at $[C]$ and
give a sufficient condition for $\Hilb^{sc} X$
to contain a generically non-reduced irreducible component of Mumford type.
\end{abstract}
\maketitle

\section{Introduction}
\label{sect:introduction}

We work over an algebraically closed field $k$ of characteristic $0$.
Given a projective scheme $X$ over $k$, 
we denote by $\Hilb^{sc} X$
the Hilbert scheme of smooth connected curves in $X$.
Mumford~\cite{Mumford} first proved that 
$\Hilb^{sc} \mathbb P^3$
contains a generically non-reduced (irreducible) component.
In \cite{Mukai-Nasu,Nasu4,Nasu5} for smooth Fano $3$-folds $X$,
$\Hilb^{sc} X$ has been studied
from the viewpoint of generalizations of Mumford's example
and more recently it has been proved in \cite{Nasu6} that
if $X$ is a prime Fano $3$-fold
then $\Hilb^{sc} X$ contains
a generically non-reduced component
whose general member is contained in a smooth hyperplane section 
$S\sim -K_X$ in $\Pic X$ ($=\mathbb Z[-K_X]$), 
i.e.~a smooth $K3$ surface $S$ in $X$.

In this paper, we study the Hilbert scheme $\Hilb^{sc} X$
for Enriques-Fano $3$-folds $X$
(see Definition~\ref{dfn:EF3})
and discuss the existence of its generically non-reduced components.
The $3$-folds $X$ in this class contain a (smooth) Enriques surface $S$
as a hyperplane section, and
were originally studied by Fano in his famous paper \cite{Fano38} and 
the study was followed in many papers, e.g,~\cite{Conte83,Conte-Murre85,
Bayle94,Sano95,Prokhorov07,Knutsen-Lopez-Munoz11}.
It is known that every Enriques-Fano $3$-fold $X$ has isolated singularities
and $-K_X$ is not a Cartier divisor but numerically 
equivalent to the hyperplane section $S \in |\mathcal O_X(1)|$.
The number $g:=(-K_X)^3/2+1$ is called 
the {\em genus} of $X$, and it is known that for every $X$ we have $g \le 17$
(cf.~\cite{Prokhorov07,Knutsen-Lopez-Munoz11}).
It follows from a general theory that on every Enriques surface $S$
there exists an effective divisor $E$ on $S$ such that
$|2E|$ is base-point-free and defines an elliptic fibration on $S$.
(Such a divisor $E$ is called a {\em half pencil} on $S$.)
Let $N_{E/X}$ denote the normal bundle of $E$ in $X$
and let $N_{E/X}(E)$ be defined by $N_{E/X}(E):=N_{E/X}\otimes_E N_{E/S}$.
The following is our main theorem.
\begin{thm}
  \label{thm:main1}
  Let $X$ be an Enriques-Fano $3$-fold of genus $g$, 
  and $S \subset X$ a smooth hyperplane section, 
  i.e.~an Enriques surface $S$ in $X$.
  If there exists a half pencil $E$ on $S$ of 
  degree $e:=(-K_X.E)\ge 2$ such that $H^1(E,N_{E/X}(E))=0$,
  then $\Hilb^{sc} X$ contains a generically non-reduced 
  component $W$ of dimension $2g+2e-1$
  whose general member $C$ satisfies:
  \begin{enumerate}
    \item $C$ is contained in an Enriques surface $S'$ in $X$, and
    \item $C$ is linearly equivalent to $-K_X\big{\vert}_{S'}+2E'$ 
    in $S'$ for some half pencil $E'$ on $S'$.
  \end{enumerate}
\end{thm}
In Example~\ref{ex:example of GNRC}
we give a few examples of Enriques-Fano $3$-folds $X$ (of genus $g=6,9,13$)
satisfying the assumption of Theorem~\ref{thm:main1}.
For these $X$, there exist a smooth Fano $3$-fold $Y$
and a $K3$ surface $M \subset Y$ 
which double cover $X$ and $S$, respectively.
We use the geometry of elliptic fibrations on $M$
to show the existence of the desired half pencil $E$ on $S$.
It might be notable that for these $X$ we have
\begin{enumerate}
  \item every general member $C$ of the component $W$
  is contained in a general hyperplane section $S$ of $X$
  (cf.~Remark~\ref{rmk:general section}), and
  \item for the smooth Fano cover $Y$ of $X$,
  there exists a generically non-reduced component $V$
  of $\Hilb^{sc} Y$ of double dimension as $\Hilb^{sc} X$
  (i.e.~$\dim V=2\dim W$ for the component $W \subset \Hilb^{sc} X$),
  whose general member is contained in a $K3$ surface $M$ in $Y$,
  but $M$ is not general in $|-K_Y|$
  (cf.~Remark~\ref{rmk:double dimension}).
\end{enumerate}
One can compare Theorem~\ref{thm:main1} with
Proposition~\ref{prop:existence of GNRC for smooth Fano},
which gives a sufficient condition for 
the Hilbert scheme of a smooth Fano $3$-fold
to have a generically non-reduced component.
Theorem~\ref{thm:main1} is obtained as an application of 
Theorem~\ref{thm:main2},
which enables us to compute the dimension of $\Hilb^{sc} X$
at $[C]$ and determines the (un)obstructedness of $C$ in $X$
for curves $C$ contained in $S$.

\begin{thm}
  \label{thm:main2}
  Let $X$ be an Enriques-Fano $3$-fold of genus $g$, 
  $S$ a smooth hyperplane section of $X$
  and $C$ a smooth connected curve on $S$
  satisfying $H^1(S,\mathcal O_S(C))=0$.
  We define a divisor $D$ on $S$ by $D:=C+K_X\big{\vert}_S$.
  \begin{enumerate}
    \item If $H^1(S,D)=0$, then $C$ is unobstructed in $X$.
    \item If there exists 
    a half-pencil $E$ on $S$ such that 
    $D \sim 2mE$ or $D \sim K_S+(2m+1)E$ for an integer $m \ge 1$,
    then we have $h^1(S,D)=m$.
    If moreover $H^1(E,N_{E/X}(E))=0$,
    then $C$ is obstructed in $X$.
    \item If $D \ge 0$, $D^2 \ge 0$
    and there exists a $(-2)$-curve $E$ on $S$
    such that $E.D=-2$ and $H^1(S,D-3E)=0$,
    then we have $h^1(S,D)=1$.
    If moreover the $\pi$-map $\pi_{E/S}(E)$
    (cf.~\eqref{map:pi-map}) for $(E,S)$ is not surjective,
    then $C$ is obstructed in $X$.
  \end{enumerate}
In (1), (2) and (3), if we assume furthermore that
$h^0(S,K_S-D)=0$ and $m=1$ (for (2)),
then $\Hilb^{sc} X$ is of dimension $g+g(C)-1$ at $[C]$,
where $g(C)$ denote the (arithmetic) genus of $C$.
\end{thm}
If $H^1(S,\mathcal O_S(C))=0$, then
the Hilbert-flag scheme $\HF^{sc} X$ of $X$ is nonsingular at $(C,S)$
of expected dimension $g+g(C)-1$
(cf.~Lemma~\ref{lem:nonsingularity of flag}).
If moreover $H^1(S,D)=0$, then the first projection 
$\HF^{sc} X \rightarrow \Hilb^{sc} X$, $(C,S)\mapsto [C]$
is smooth at $(C,S)$, and thus Theorem~\ref{thm:main2}~(1) follows
from a property of smooth morphisms.
We partially prove that $C$ is obstructed in $X$ if $H^1(S,D)\ne 0$
by using half pencils and $(-2)$-curves on $S$
together with a result in \cite{Nasu5}.
See \cite{Nasu5} for a result on the (un)obstructedness of curves 
lying on a $K3$ surface in a smooth Fano $3$-fold.

The organization of this paper is as follows.
In \S\ref{subsec:enriques} and \S\ref{subsec:EF3}
we recall some properties of Enriques surfaces and Enriques-Fano $3$-folds,
respectively.
In \S\ref{subsec:flag and primary obstruction}
we recall some known results on Hilbert-flag schemes
and obstructions to lifting first order deformations
of curves on a $3$-fold to second order deformations
(i.e.~priamary obstructions).
These results will be used in 
\S\ref{sect:obstructions} and \S\ref{sect:non-reduced}
to prove Theorems~\ref{thm:main2} and \ref{thm:main1}, respectively.

\smallskip

\paragraph{\bf Acknowledgments}
I would like to thank Prof.~Hiromichi Takagi
for his comment, which motivated me to research the topic of this paper.
I would like to thank Prof.~Shigeru Mukai
for letting me know examples of Enriques-Fano $3$-folds.
This paper was written during my stay as a visiting researcher
at the department of mathematics at the University of Oslo (UiO), Norway.
I thank UiO for providing the facilities.
I thank Prof.~Kristian Ranestad, Prof.~John Christian Ottem and
Prof.~Jan Oddvar Kleppe for helpful and inspiring discussions during the stay.
Last but not least, I thank the referee for giving helpful comments improving the readability and quality of this paper. 
This work was supported in part by JSPS KAKENHI Grant Numbers 
JP17K05210 and JP20K03541.

\section{Preliminaries}
\label{sect:preliminaries}

\subsection{Enriques surfaces}
\label{subsec:enriques}

In this section, we recall some properties of Enriques surfaces.
We refer to \cite{Dolgachev16} and \cite{BHPV}
for proofs and more general theories on Enriques surfaces.
A smooth projective surface $S$
is called an {\em Enriques surface}
if $H^i(S,\mathcal O_S)=0$ for $i=1,2$ and $2K_S\sim 0$.
Every Enriques surface $S$ is isomorphic to
the quotient $M/\theta$ of a smooth $K3$ surface $M$
by a fixed-point-free involution $\theta$ of $M$.
Here $M$ is the canonical cover of $S$ and
called the {\em $K3$ cover} of $S$.
It is well known that $S$ admits an elliptic fibration
$\psi: S \longrightarrow \mathbb P^1$,
whose general fiber is a smooth curve of arithmetic genus $1$
(cf.~e.g.~\cite[VIII,~\S17]{BHPV}).
For each $\psi$
there exist exactly two multiple fibers $2E$ and $2E'$ of $\psi$,
i.e.~the {\em double fibers} of $\psi$,
and we have $K_S\simeq \mathcal O_S(E-E')\simeq \mathcal O_S(E'-E)$.
Such divisors $E$ and $E'$ are called {\em half pencils} on $S$.
Every nef and primitive divisor $E\ge 0$ on $S$
with self-intersection number $E^2=0$ is a half pencil on $S$.

Let $S$ be an Enriques surface and 
$\pi: M \rightarrow S$ its $K3$ cover, 
i.e.~there exists a fixed-point-free involution $\theta$ of $M$
such that $S \simeq M/\theta$.
Let $\varphi: M\rightarrow \mathbb P^1$ 
be an elliptic fibration on $M$. 
Then there exists a non-constant rational function $u \in k(M)$,
which is defined by $u=\varphi^*t$ and $k(\mathbb P^1)=k(t)$.
Here $u$ is called {\em the elliptic parameter} 
of $\varphi$ 
and unique only up to the linear fractional transformations (of $t$)
(cf.~\cite{Kuwata-Shioda08}).

\begin{lem}
  \label{lem:elliptic fibration}
  Let $u$ be the elliptic parameter of $\varphi: M \rightarrow \mathbb P^1$
  and $F$ the fiber of $\varphi$ defined by $u$,
  i.e.~$F=\div_0(u)$.
  If $u$ is $\theta$-anti-invariant, i.e.~$\theta^* u=-u$,
  then the image $E=\pi(F)$ is a half pencil on $S$
  and the pull-back $\pi^*E$ coincides with $F$.
\end{lem}

\Proof Let $\pi'$ be the map on $\mathbb P^1$ defined by $t \mapsto t^2$.
Since $u^2$ is invariant under $\theta$, the composition $\pi'\circ \varphi$
factors through $\pi$ and there exists a commutative diagram
$$
\begin{CD}
  M @>{\varphi}>> \mathbb P^1 \\
  @V{\pi}VV @V{\pi'}VV \\
  S @>{\psi}>> \mathbb P^1,
\end{CD}
$$
where $\psi$ defines an elliptic fibration on $S$.
Since $\pi'$ is ramified at $(t=0)$ (and $(t=\infty)$) on $\mathbb P^1$, 
$\psi$ has the double fibers $2E$ at $(t=0)$. 
By commutativity, we have proved the lemma.
\qed

\medskip

Let $E$ be a half pencil on $S$.
Note that the normal bundle $N_{E/S}$ ($\simeq \mathcal O_E(E)$) 
of $E$ in $S$ is a $2$-torsion in $\Pic E$.
Here $\mathcal O_E(E)$ has no sections, but 
$\mathcal O_E(2E)$ has a unique nonzero section up to constants.
We also note that for every integer $k \ge 0$,
the restriction map $H^0(S,kE) \rightarrow H^0(E,kE)$
is surjective.
It follows from an exact sequence
$0 \rightarrow \mathcal O_S((k-1)E)
\rightarrow \mathcal O_S(kE)
\rightarrow \mathcal O_E(kE)
\rightarrow 0$ that
the surjectivity is equivalent to the injectivity of
the natural induced map 
$H^1(S,(k-1)E) \rightarrow H^1(S,kE)$.
We need the following lemma for our proof of Theorem~\ref{thm:main2}.

\begin{lem}
  \label{lem:KL.etc}
  Let $D$ be an effective divisor on $S$ with $D^2 \ge 0$.
  \begin{enumerate}
    \item $H^1(S,D)\ne 0$ if and only if
    \begin{enumerate}
      \item $D \sim mE$ or $D \sim K_S+(m+1)E$ 
      for a half pencil $E$ on $S$ and an integer $m \ge 2$
      (then $h^1(S,D)=\lfloor m/2 \rfloor$), or
      \item $D.\Delta \le -2$
      for some divisor $\Delta \ge 0$ on $S$ with $\Delta^2 =-2$.
    \end{enumerate}
    \item If there exists $(-2)$-curve $E$ on $S$ such that $E.D=-2$
    and $H^1(S,D-3E)=0$ then we have $h^1(S,D)=1$ and $H^1(S,D-E)=0$.
  \end{enumerate}
\end{lem}

\Proof
(1) is a special case of \cite{Knutsen-Lopez}. (2) follows from 
\cite[Claim~4.1]{Nasu5}, whose proof works also for
divisors on Enriques surfaces.
\qed

\subsection{Enriques-Fano $3$-folds}
\label{subsec:EF3}

In this section, we collect some known results on 
Enriques-Fano $3$-folds and prepare a lemma
on the elliptic fibrations on their hyperplane sections
(cf.~Lemma~\ref{lem:elliptic fibration on K3}).
This lemma will be used in \S\ref{sect:non-reduced}
to show the existence of a generically non-reduced component
of the Hilbert scheme of some Enriques-Fano $3$-folds.

\begin{dfn}
  \label{dfn:EF3}
  A normal projective $3$-dimensional variety 
  $X \subset \mathbb P^N$ is called 
  {\em Enriques-Fano} if it contains
  an Enriques surface $S$ as a hyperplane section 
  and $X$ is not a cone over $S$.
\end{dfn}
An equivalent definition is to assume that 
a general hyperplane section is a smooth Enriques surface.
In this section, we mainly consider the Enriques-Fano $3$-folds $X$ with
only terminal cyclic quotient singularities.
There is a classification of such $X$ due to
Bayle~\cite{Bayle94} and Sano~\cite{Sano95}.
We summarize the properties of $X$:
\begin{enumerate}
  \item $-2K_X \sim 2H$ in $\Pic X$ for $H \in |\mathcal O_X(1)|$. 
  \item The canonical cover
  $$
  Y=\Spec_X (\mathcal O_X \oplus \mathcal O_X(K_X+H))
  \overset{\pi}{\longrightarrow} X
  $$
  of $X$ is a smooth Fano $3$-fold,
  and $Y$ is isomorphic to one of the $3$-folds
  in Table \ref{table:EF3}.
  \item The covering transformation $\theta$ of $\pi$
  is an involution of $Y$
  and $\theta$ fixes just $8$ points on $Y$.
  \item The fixed points on $Y$ give rise to
  the singularity on $X \simeq Y/\theta$ of type $\frac 12(1,1,1)$.
\end{enumerate}

\begin{table}[h]
  \caption{Enriques-Fano $3$-folds $X$}
  \begin{center}
    \begin{minipage}{14.5cm}
      \begin{tabular}{|c|l|c|}
	\hline
	No. & canonical cover $Y$
	$(\overset{2:1}{\longrightarrow} X)$ & $g$\\
      \hline
      1 & a complete intersection
      $(2) \cap (4) \subset \mathbb P(1,1,1,1,1,2)$  & $2$ \\
      \hline
      2 & a complete intersection 
      $(2) \cap (2) \cap (2) \subset \mathbb P^6$ & $3$ \\
      \hline
      3 & the blow-up $\Bl_{\gamma} V_2$ of $V_2$\footnote{
	In this table, 
	$V_n$ denotes a del Pezzo $3$-fold of degree $n$.}
      with a center an elliptic curve $\gamma\footnote{
	$\gamma$ is a complete intersection $H_1\cap H_2 \subset V_2$
      where $H_i \in |(-1/2)K_{V_2}|$ for $i=1,2$.} \subset V_2$ & $3$\\
      \hline
      4 & $\mathbb P^1 \times S_2$\footnote{
	In this table, $S_n$ denotes a del Pezzo surface of degree $n$.} 
      & $4$ \\
      \hline
      5 & a double cover of $\mathbb P^1\times \mathbb P^1 \times \mathbb P^1$
      branched along a divisor $R\sim (2,2,2)$. & $4$ \\
      \hline
      6 & a double cover of $(1,1) \subset \mathbb P^2\times \mathbb P^2$
      \footnote{
	A hypersurface in $\mathbb P^2\times \mathbb P^2$
	of multidegree $(1,1)$ is isomorphic to
	$\mathbb P_{\mathbb P^2}(T_{\mathbb P^2})$.}
      branched along $R \in |-K|$. & $4$ \\
      \hline
      7 & the blow-up $\Bl_{\gamma} V_4$ of $V_4$
      with a center an elliptic curve
      $\gamma\footnote{
	$\gamma$ is a complete intersection 
	$(1)\cap (1)\cap (2)\cap (2) \subset \mathbb P^5$} 
      \subset \mathbb P^5$
      & $5$\\
      \hline
      8 & a weighted hypersurface $(4)\subset \mathbb P(1,1,1,1,2)$ 
      (i.e.~$V_2$) & $5$ \\
      \hline
      9 & a complete intersection $(1,1)\cap (1,1)\cap (1,1) 
      \subset \mathbb P^3\times\mathbb P^3$ & $6$ \\
      \hline
      10 & $\mathbb P^1\times S_4$ & $6$ \\
      \hline
      11 & a hypersurface in 
      $\mathbb P^1\times \mathbb P^1\times \mathbb P^1\times \mathbb P^1$
      of multidegree $(1,1,1,1)$
      & $7$ \\
      \hline
      12 & a complete intersection 
      $(2)\cap (2)\subset \mathbb P^5$ (i.e.~$V_4$) & $9$ \\
      \hline
      13 & $\mathbb P^1\times S_6$ & $10$ \\
      \hline
      14 & $\mathbb P^1\times \mathbb P^1 \times \mathbb P^1$ & $13$\\
      \hline
    \end{tabular}
  \end{minipage}
\end{center}
  \label{table:EF3}
\end{table}

Minagawa~\cite{Minagawa99} showed that
every Enriques-Fano $3$-fold $X_0$ with at most terminal singularities
admits a $\mathbb Q$-smoothing, which assures us that
$X_0$ is obtained as a flat specialization of Enriques-Fano $3$-folds $X_t$
($t \ne 0$) with only terminal cyclic quotient singularities.
This fact implies that we have $g \le 13$ for any $X_0$.
\begin{rmk}
  \label{rmk:classification of EF3}
  Prokhorov~\cite{Prokhorov07} and 
  Knutsen-Lopez-Mu\~{n}oz~\cite{Knutsen-Lopez-Munoz11} proved that $g\le 17$
  for any Enriques-Fano $3$-fold $X$.
  As far as we know, the problem of classifying
  all Enriques-Fano $3$-folds $X$ with non-terminal singularities is still open.
\end{rmk}

In what follows, we recall some well known examples
of Enriques-Fano $3$-folds $X_g$ of genus $g=6,9$ and $13$.
Here $X_g$ has the Picard rank $1,2$ and $3$ 
for $g=9,6$ and $13$, respectively.
\begin{ex}
  \label{ex:EF3}
  In the following example, $Y$ is a smooth Fano $3$-fold
  and there exists an involution $\theta$ of $Y$ fixing just 
  $8$ points on $Y$.
  There exists a smooth $K3$ surface $M$ in $Y$
  on which $\theta_M:=\theta\big{\vert}_M$ acts without fixed points.
  Thereby the quotient $S:=M/\theta_M$ is an Enriques surface
  and $X:=Y/\theta$ is an Enriques-Fano $3$-fold of genus 
  $g=(-K_X)^3/2+1=(-K_Y)^3/4+1$.

  \begin{enumerate}
    \item (No.14 in Table~\ref{table:EF3})
    Let $Y :=\mathbb P^1 \times \mathbb P^1 \times \mathbb P^1$.
    We define an involution $\theta$ of $Y$ by
    \begin{equation}
      \label{map:involution (-1,-1,-1)}
      (x_0:x_1)\times(y_0:y_1)\times(z_0:z_1)
      \longmapsto
      (x_0:-x_1)\times(y_0:-y_1)\times(z_0:-z_1).
    \end{equation}
    Here we say that $\theta$ is {\em of type} $(-1,-1,-1)$.
    Then $\theta$ fixes just $8$ coordinate points on $Y$.
    There exist exactly $14$ $\theta$-invariant monomials,
    which correspond to the $14$ ($=\bullet \times 8 +\circ \times 6$)
    vertices in Figure~\ref{fig:invariant(2,2,2)}. 
    Then these monomials span the linear subsystem $\Lambda$ of
    $|-K_Y|=|\mathcal O_Y(2,2,2)|\simeq \mathbb P^{26}$
    of dimension $13$.
    \begin{figure}[h]
      \begin{center}
	\subfigure[invariant $(2,2,2)$-forms]{
	  \raisebox{-1cm}{\includegraphics[clip,scale=0.7]{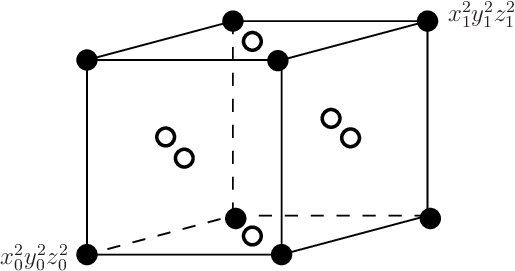}
	  \label{fig:invariant(2,2,2)}
	}}
	\subfigure[(anti-)invariant $(1,1,1)$-forms]{
	\raisebox{-1cm}{\includegraphics[clip,scale=0.6]{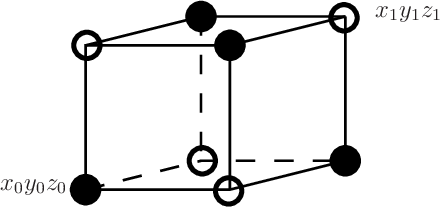}
	  \label{fig:(anti)invariant(1,1,1)}
	  }}
      \end{center}
      \caption{(anti-)invariant monomials}
    \end{figure}
    There exists a smooth member $M \in \Lambda$
    not passing through the $8$ fixed points.
    Thus $X_{13}:=Y/\theta$ is an Enriques-Fano $3$-fold of genus $13$,
    and $X_{13}$ is embedded into 
    $\mathbb P^{13}$ by the linear system
    $\Lambda\simeq |\mathcal O_X(2,2,2)|\simeq \mathbb P^{13}$.

    \item (No.12 in Table~\ref{table:EF3})
    Let $\mathbb P^5$ be the projective $5$-space
    and $x_0,\dots,x_5$ its homogeneous coordinates.
    Let $Y\subset \mathbb P^5$ be a smooth complete intersection 
    of two hyperquadrics, whose defining equations are of the forms
    \begin{equation}
      \label{eqn:quadrics}
      q_i(x_0,x_1,x_2)+q'_i(x_3,x_4,x_5)
      \qquad (i=0,1).
    \end{equation}
    We define an involution $\theta'$ of $\mathbb P^5$ by
    $$
    (x_0,x_1,x_2,x_3,x_4,x_5) \longmapsto 
    (x_0,x_1,x_2,-x_3,-x_4,-x_5).
    $$
    Then $\theta'$ defines an involution $\theta$ of $Y$ by restriction.
    The the fixed locus $\Fix(\theta')$ of $\theta'$
    is equal to $(x_0=x_1=x_2=0)\bigcup (x_3=x_4=x_5=0)$.
    Thereby $\theta$ fixes just $8$ points on $Y$.
    We consider the third quadratic form
    $q_2(x_0,x_1,x_2)+q_2'(x_3,x_4,x_5)$
    of the same type and 
    the hyperquadric $Q_2$ in $\mathbb P^5$ defined by it.
    Then the intersection $M:=Y \cap Q_2$ is a $K3$ surface,
    and moreover we can take $q_2$ and $q_2'$ so that 
    $M \cap \Fix(\theta)=\emptyset$.
    Thus $X_9:=Y/\theta$ is an Enriques-Fano $3$-fold of genus $9$.
    Then the linear system $\Lambda$ $(\simeq \mathbb P^9)$
    of $\theta$-invariant quadratic forms on $Y$ of type \eqref{eqn:quadrics}
    defines an embedding of $X_9$ into $\mathbb P^9$.

    \item (No.9 in Table~\ref{table:EF3})
    Let $Y$ be a smooth complete intersection
    of three hypersurfaces $Q_i$ ($i=0,1,2$)
    in $\mathbb P^3\times \mathbb P^3$ of bidegree $(1,1)$.
    Suppose that for each $i$,
    $Q_i$ is defined by a symmetric bilinear form
    $q_i(\mathbf x,\mathbf y)$ on $\mathbb P^3\times \mathbb P^3$.
    We consider the diagonal action
    $(\mathbf x,\mathbf y) \longmapsto (\mathbf y,\mathbf x)$
    on $\mathbb P^3\times \mathbb P^3$
    and define an involution $\theta$ on $Y$ by its restriction.
    Then $\theta$ fixes just $8$ diagonal points on $Y$.
    We consider the fourth symmetric bilinear form 
    $q_3(\mathbf x,\mathbf y)$ and the hypersurface 
    $Q_3 \subset \mathbb P^3\times \mathbb P^3$ defined by it.
    If $q_3$ is general, then the intersection $M:=Y \cap Q_3$ 
    is a $K3$ surface in $Y$,
    on which $\theta\big{\vert}_M$ acts without fixed points.
    Thus $X_6:=Y/\theta$ is an Enriques-Fano $3$-fold of genus $6$.
    The linear system $\Lambda$ ($\simeq \mathbb P^6$)
    of symmetric bilinear forms on $Y$
    defines an embedding of $X_6$ into $\mathbb P^6$.
  \end{enumerate}
\end{ex}

Let $X$ be an Enriques-Fano $3$-fold
with at most terminal cyclic quotient singularities.
Then the canonical (double) cover $\pi: Y \rightarrow X$ of $X$ 
is a smooth Fano $3$-fold and there exists
an involution $\theta$ of $Y$ such that $X \simeq Y/\theta$.
(Here and later we call $Y$ the {\em smooth Fano cover} of $X$.)
If a $K3$ surface $M \subset Y$ is invariant under $\theta$
and an elliptic fibration $\varphi:M \rightarrow \mathbb P^1$ 
has a $\theta$-anti-invariant elliptic parameter $u \in k(M)$,
then by Lemma~\ref{lem:elliptic fibration},
$\varphi$ induces an elliptic fibration on 
the Enriques quotient $S \subset X$ of $M$,
and the image $E=\pi(F)$ of $F=(u=0)$ becomes a half-pencil on $S$.

\begin{lem}
  \label{lem:elliptic fibration on K3}
  If $X$ is either $X_{13},X_9$ or $X_6$ in Example~\ref{ex:EF3},
  then there exist a $\theta$-invariant $K3$ surface $M$ in 
  the smooth Fano cover $Y$ of $X$
  and an elliptic fibration 
  $\varphi: M \rightarrow \mathbb P^1$
  with a $\theta$-anti-invariant elliptic parameter $u$
  such that 
  $$
  H^1(F,N_{F/Y})=0
  $$
  for the invariant elliptic curve $F=(u=0)$ on $M$.
  More explicitly, $F$ is described as follows:
  \begin{enumerate}
    \item $F$ is a complete intersection of two hypersurfaces
    of tridegree $(1,1,1)$ in 
    $Y=\mathbb P^1\times\mathbb P^1\times\mathbb P^1$
    if $X=X_{13}$,
    \item $F$ is a linear section $(1)\cap (1)\cap Y$ of 
    $Y=(2)\cap (2) \subset \mathbb P^5$ if $X=X_9$, and
    \item $F$ is a complete intersection of two hypersurfaces
    of bidegree $(1,0)$ and $(0,1)$ with 
    $Y=(1,1)^3 \subset \mathbb P^3 \times \mathbb P^3$ if $X=X_6$.
  \end{enumerate}
\end{lem}

\Proof
The lemma is proved in case by case.

\noindent
(1)\quad
Let $X=X_{13}$. 
We denote by $R^+$ (resp.~$R^-$)
the space of invariant (resp.~anti-invariant) $(1,1,1)$-forms on $Y$.
Then $R^+$ and $R^-$ are spanned by
the four monomials $x_iy_jz_k$ corresponding to
the vertices $\bullet$ and $\circ$
in Figure~\ref{fig:(anti)invariant(1,1,1)}, respectively.
We take four general $(1,1,1)$-forms 
$l^+,m^+ \in R^+$ and $l^-,m^- \in R^-$
and define a smooth $K3$ surface $M$ in $Y$
by a $(2,2,2)$-form $f=l^+m^+ - l^-m^-$ on $Y$.
Then since $f$ is $\theta$-invariant, so is $M$.
We consider a rational function $u=l^+/m^-$ ($=l^- / m^+$) on $M$
and the elliptic fibration $\varphi: M\rightarrow \mathbb P^1$
with the elliptic parameter $u$.
Then $u$ is clearly $\theta$-anti-invariant.
The fiber $F$ at $(t=0) \in \mathbb P^1$ is
defined by $l^+=l^-=0$ in $Y$, and hence $F$ is a complete intersection in $Y$.
Since $N_{F/Y}\simeq \mathcal O_F(1,1,1)^{\oplus 2}$ and 
$\mathcal O_F(1,1,1)$ is ample, we have $H^1(F,N_{F/Y})=0$.

\noindent
(2)\quad
Let $X=X_9$ and let $M$ be the $K3$ surface in Example~\ref{ex:EF3}(2).
Then $M$ is defined in $\mathbb P^5$ by
three quadratic polynomials $q_i(\mathbf x) + q_i'(\mathbf x')$ 
($i=0,1,2$), where $\mathbf x=(x_0,x_1,x_2)$ and $\mathbf x'=(x_3,x_4,x_5)$.
Then for each $i$,
there exist two $3\times 3$ symmetric matrices $A_i$ and $A_i'$
corresponding to $q_i$ and $q_i'$, respectively.
We see that there is a one-to-one correspondence between 
the set of elliptic fibrations on $M$ (on $S=M/\theta$)
and the $9$ points
$\lambda=(\lambda_0,\lambda_1,\lambda_2)$ in $\mathbb P^2$ 
defined by 
$$
\det (\lambda_0 A_0+\lambda_1 A_1+\lambda_2 A_2)
=\det (\lambda_0 A'_0+\lambda_1 A'_1+\lambda_2 A'_2)
=0. 
$$
In fact, if $\lambda$ satisfies this equation,
then there exist four linear forms $l,m,l',m'$ on $\mathbb P^2$
such that
$\sum_{i=0}^2 \lambda_i (q_i(\mathbf x)+q_i'(\mathbf x'))
=l(\mathbf x)m(\mathbf x)-l'(\mathbf x')m'(\mathbf x')$.
Then as in (1), there exists an elliptic fibration 
$\varphi: M \rightarrow \mathbb P^1$ on $M$ defined by 
the elliptic parameter
$u=l(\mathbf x)/m'(\mathbf x')$ ($=l'(\mathbf x')/m(\mathbf x)$),
which is $\theta$-anti-invariant.
Since the fiber $F=(u=0)$ of $\varphi$ is 
defined by $l(\mathbf x)=l'(\mathbf x')=0$ in $Y$,
$F$ is a linear section of $Y$.
Thus $H^1(F,N_{F/Y})=H^1(F,\mathcal O_F(1)^{\oplus 2})=0$.

\noindent
(3)\quad
Let $X=X_6$ 
and let $M$ be the $K3$ surface in Example~\ref{ex:EF3}(3).
We recall that $M$ is a complete intersection
of four hypersurfaces $Q_i \subset \mathbb P^3 \times \mathbb P^3$ 
($i=0,1,2,3$) of bidegree $(1,1)$, each of which
is defined by a symmetric bilinear form $q_i$ 
on $\mathbb P^3 \times \mathbb P^3$.
For each $i$, let $A_i$ be the $4\times 4$ symmetric matrix
corresponding to $q_i$.
Then the elliptic fibrations on $M$ (on $S=M/\theta$)
are in one-to-one correspondence with 
the points $\lambda=(\lambda_0,\lambda_1,\lambda_2,\lambda_3)$ 
in $\mathbb P^3$ satisfying
$$
\rank (A(\lambda)) \le 2,
\qquad
\text{where $A(\lambda):=\lambda_0 A_0 +\lambda_1 A_1 +\lambda_2 A_2 + \lambda_3 A_3$}.
$$
These points correspond to the 10 nodes of 
a quartic surface in $\mathbb P^3$ defined by $\det (A(\lambda))$, 
which is well known as a Cayley's quartic symmetroid (cf.~\cite{Cayley69}).
For each such $\lambda$,
there exist two linear forms $l$ and $m$ on $\mathbb P^3$ such that
${}^t \mathbf x A(\lambda) \mathbf y=
l(\mathbf x)l(\mathbf y)+m(\mathbf x)m(\mathbf y)$.
By changing the coordinates of $\mathbb P^3$, we may assume that
${}^t \mathbf x A(\lambda) \mathbf y=x_0y_0+x_1y_1$.
We define a rational function $u \in k(M)$
by $u=(x_1-x_0\sqrt{-1})/(x_1+x_0\sqrt{-1})$
and take it as an elliptic parameter on $M$.
It is easy to see that $u$ is anti-invariant under
the involution $\theta(\mathbf x,\mathbf y)=(\mathbf y,\mathbf x)$.
Moreover every fiber of the elliptic fibration defined by $u$
is a complete intersection $(1,0) \cap (0,1) \cap Y$ in $Y$.
Thus we have $N_{F/Y}
\simeq \mathcal O_F(1,0)\oplus \mathcal O_F(0,1)
\simeq \pi_1^*\mathcal O_{\mathbb P^3}(1)\oplus 
\pi_2^*\mathcal O_{\mathbb P^3}(1)$,
where $\pi_i: F \rightarrow \mathbb P^3$ is the restriction to $F$
of the $i$-th projection 
$\mathbb P^3\times \mathbb P^3 \rightarrow \mathbb P^3$ 
($i=1,2$). This implies that $H^1(F,N_{F/Y})=0$.
\qed

\begin{rmk}
  \label{rmk:general section}
  We note that in Lemma~\ref{lem:elliptic fibration on K3},
  the $K3$ surface $M$ is required to
  contain a pencil of elliptic curves $F$ in $Y$
  and hence not a general member of $|-K_Y|$.
  However its image $S=\pi(M)$ is a general hyperplane section of
  $X_g \subset \mathbb P^g$.
  It is rather easy to see this for $g=6,9$.
  For $g=13$ we consider a linear map
  $$
  \Phi:
  \Sym^2(R^+) \oplus \Sym^2(R^-) 
  \longrightarrow H^0(X,\mathcal O_X(2,2,2))\simeq k^{14},
  $$
  where $R^+$ (resp.~$R^-$) is the $4$-dimensional vector space
  of (resp.~anti-)invariant $(1,1,1)$-forms on $Y$ 
  (cf.~Figure~\ref{fig:(anti)invariant(1,1,1)}), and
  $\Sym^2(V)$ denotes the symmetric square of a $k$-vector space $V$.
  Then we see that the kernel of $\Phi$ is of dimension 
  $6$ and hence $\Phi$ is surjective.
\end{rmk}

\begin{rmk}
  \label{rmk:-K-normal}
  In Lemma~\ref{lem:elliptic fibration on K3},
  $F$ is a complete intersection in $Y$.
  This implies that
  $H^1(Y,\mathcal I_{F/Y}\otimes_Y \mathcal O_Y(-K_Y))=0$.
  Then there exists a short exact sequence
  $$
  \begin{CD}
    0 @>>> \mathcal O_Y
    @>>> \mathcal I_{F/Y}\otimes_Y \mathcal O_Y(-K_Y)
    @>>> N_{M/Y}(-F)
    @>>> 0
  \end{CD}
  $$
  on $Y$ and it follows from this sequence that $H^1(M,N_{M/Y}(-F))=0$.
\end{rmk}

\subsection{Hilbert-flag scheme and Primary obstruction}
\label{subsec:flag and primary obstruction}

In this section, 
we recall some properties of Hilbert-flag schemes
(cf.~\cite{Kleppe81, Sernesi})
and a result in \cite{Nasu5} on primary obstructions
to deforming curves on a $3$-fold.
In this paper, Hilbert-flag schemes play an important role.
Given a projective scheme $X$, we denote by $\HF X$ 
the {\em Hilbert-flag scheme}, that is the projective scheme
parametrizing all pairs $(C,S)$ of closed subschemes $C$ and 
$S$ of $X$ satisfying $C \subset S$. 
For a pair $(C,S)$, its normal sheaf $N_{(C,S)/X}$ in $X$ 
is defined by the Cartesian diagram
$$
\vcenter{
  \xymatrix{
    N_{(C,S)/X} \ar[d]_{\pi_1} \ar[r]^{\pi_2} \ar@{}[dr]|\square & 
    N_{S/X} \ar[d]_{|_C} \\
    N_{C/X} \ar[r]^{\pi_{C/S}} & N_{S/X}\big{\vert}_C, \\
}}
$$ 
where $|_C$ and $\pi_{C/S}$ are the restriction and 
the projection, respectively.
Suppose that the two embeddings $C \hookrightarrow S$ and 
$S \hookrightarrow X$ are regular (cf.~\cite{Sernesi}).
Then $H^0(X,N_{(C,S)/X})$ and $H^1(X,N_{(C,S)/X})$ respectively 
represent the tangent space and the obstruction space
of $\HF X$ at $(C,S)$. It follows from a general theory 
on Hilbert schemes that
\begin{equation}
  \label{ineq:dimension of flag}
  h^0(X,N_{(C,S)/X})-h^1(X,N_{(C,S)/X}) \le \dim_{(C,S)} \HF X
  \le h^0(X,N_{(C,S)/X}).
\end{equation}
Moreover, there exist two fundamental exact sequences
\begin{equation}
  \label{ses:normal1}
  \begin{CD}
      0 @>>> \mathcal I_{C/S}\otimes_S N_{S/X} 
      @>>> N_{(C,S)/X} @>{\pi_1}>> N_{C/X} @>>> 0
  \end{CD}
\end{equation}
and
\begin{equation}
  \label{ses:normal2}
  \begin{CD}
  0 @>>> N_{C/S} @>>> N_{(C,S)/X} @>{\pi_2}>> N_{S/X} @>>> 0
  \end{CD}
\end{equation}
of sheaves on $X$ associated to the two projections
$pr_i: \HF X \rightarrow \Hilb X$ ($i=1,2$),
which send $(C,S)$ to $[C]$ for $i=1$ and $[S]$ for $i=2$, respectively.
Each of the two exact sequences induces 
a long exact sequence of cohomology groups,
which contains the tangent map and the map on obstruction spaces 
of Hilbert(-flag) schemes (cf.~\cite{Nasu6}).

Now let $X$ be a (possibly singular) projective $3$-fold,
$S$ a smooth surface contained in the smooth locus of $X$
and $C$ a smooth connected curve on $S$.
We consider the subscheme of $\HF X$ parametrising 
pairs $(C,S)$ of a smooth connected curve $C$ and a surface $S$ in $X$
and denote it by $\HF^{sc} X$.
Here and later, given a sheaf $\mathcal F$ on $X$,
we denote by $\chi(X,\mathcal F)$ 
the Euler-Poincare characteristic of $\mathcal F$.
\begin{lem}
  \label{lem:nonsingularity of flag}
  Suppose that $-K_X$ is ample and $p_g(S)=0$.
  \begin{enumerate}
    \item If $H^1(S,\mathcal O_S(C))=0$ then
    $H^i(X,N_{(C,S)/X})=0$ for all $i>0$,
    which implies that $\HF^{sc} X$ is nonsingular at $(C,S)$ of
    expected dimension $\chi(X,N_{(C,S)/X})$.
    \item If $H^1(S,\mathcal O_S(C))=H^1(S,\mathcal O_S(C+K_X\big{\vert}_S))=0$,
    then $C$ is unobstructed in $X$.
  \end{enumerate}
\end{lem}

\Proof 
(1)\quad 
There exists an exact sequence
$0 \rightarrow \mathcal O_S \rightarrow \mathcal O_S(C)
\rightarrow N_{C/S} \rightarrow 0$.
Since $H^2(S,\mathcal O_S)=0$, we have $H^i(C,N_{C/S})=0$ for all $i>0$.
By adjunction, we have $N_{S/X}\simeq-K_X\big{\vert}_S+K_S$ and hence
$H^i(S,N_{S/X})=0$ for all $i>0$ by assumption.
Thus we obtain the first conclusion of the lemma
by \eqref{ses:normal2} and \eqref{ineq:dimension of flag}.

(2)\quad
By Serre duality, we see that 
$$
H^1(S,N_{S/X}(-C))
\simeq H^1(S,-K_X\big{\vert}_S+K_S-C)
\simeq H^1(S,C+K_X\big{\vert}_S)^{\vee}=0.
$$ 
Then the first projection $pr_1$ is smooth at $(C,S)$
by \cite[Lemma~A10]{Kleppe87} (cf.~\cite[\S2.2]{Nasu6}).
Thus as a consequence of (1) we have proved the lemma.
\qed

\begin{rmk}
  \label{rmk:expected dimension}
  If $X$ is an Enriques-Fano $3$-fold and $S$ is its hyperplane section, 
  then by using \eqref{ses:normal2} 
  the number $\chi(X,N_{(C,S)/X})$ is computed as
  $$
  \chi(S,N_{S/X})+\chi(C,N_{C/S})
  =(-K_X)^3/2+1+C^2/2
  =g+g(C)-1,
  $$
  where $g$ is the genus $g=(-K_X)^3/2+1$ of $X$.
\end{rmk}

Next we recall primary obstructions to deforming curves on a $3$-fold.
Let $\alpha$ be a global section of $N_{C/X}$, i.e.~a first order 
(infinitesimal) deformation $\tilde C$ of $C$ in $X$.
We note that $C$ is a locally complete intersection in $X$.
Then the primary obstruction $\ob(\alpha)$,
i.e.~the obstruction to extend $\tilde C$ 
to a deformation ${\tilde {\tilde C}}$ of $C$ over $k[t]/(t^3)$,
is contained in $H^1(C,N_{C/X})$
and expressed as a cup product of cohomology classes on $C$
(cf.~e.g.~\cite[Theorem~2.1]{Nasu5}).
If $\ob(\alpha)\ne 0$ for some $\alpha \in H^0(C,N_{C/X})$,
then $C$ is obstructed in $X$.
In \cite{Mukai-Nasu, Nasu5} a sufficient condition for $\ob(\alpha)\ne 0$
was given under the presence of an intermediate smooth surface $S$
satisfying $C \subset S \subset X$.
Let $\pi_{C/S}:N_{C/X}\longrightarrow N_{S/X}\big{\vert}_C$
be the natural projection of normal bundles,
and $\pi_{C/S}(\alpha) \in H^0(S,N_{S/X}\big{\vert}_C)$ 
the image of $\alpha$ by the projection 
(i.e.~the {\em exterior component} of $\alpha$).
We suppose furthermore that 
$\pi_{C/S}(\alpha)$ lifts to 
a global section $\beta \in H^0(S,N_{S/X}(E))$
for some effective divisor $E\ge 0$ on $S$, i.e.~we have
\begin{equation}
\label{eqn:lifting}
r(\pi_{C/S}(\alpha),E)
=\beta\big{\vert}_C
\qquad
\mbox{in}
\qquad
H^0(C,N_{S/X}(E)\big{\vert}_C).
\end{equation}
Here and later, for a sheaf $\mathcal F$ on $S$ and 
a cohomology class $*$ in $H^i(S,\mathcal F)$,
we denote by $r(*,E)$ the image of $*$ by the natural map
$H^i(S,\mathcal F) \rightarrow H^i(S,\mathcal F(E))$
(and we use similar notation for $C$).
The rational section $\beta$ of $N_{S/X}$ admitting a pole along $E$
is called an {\em infinitesimal deformation of $S$ in $X$ with pole} (along $E$).
Tensoring the projection $\pi_{E/S}$ with $N_{E/S}\simeq \mathcal O_E(E)$ 
and taking the map induced on the space of global sections,
we define the {\em $\pi$-map}
\begin{equation}
  \label{map:pi-map}
  \pi_{E/S}(E): H^0(E,N_{E/X}(E)) 
  \longrightarrow H^0(E,N_{S/X}(E)\big{\vert}_E)
\end{equation}
for $(E,S)$.
The following theorem is crucial to our proof of 
Theorems~\ref{thm:main1} and \ref{thm:main2}.

\begin{thm}[{cf.~\cite[Theorem~1.1 and Corollary~3.2]{Nasu5}}]
  \label{thm:obstrcutedness theorem}
  The primary obstruction $\ob(\alpha)$ is nonzero if
  \begin{enumerate}
    \renewcommand{\theenumi}{{\roman{enumi}}}
    \renewcommand{\labelenumi}{{\rm (\theenumi)}}
    \item the natural map $H^1(S,kE)\longrightarrow H^1(S,(k+1)E)$
    is injective for all integer $k\ge 0$,
    \label{item:injective}
    \item the restriction map
    $H^0(S,\Delta)\longrightarrow H^0(E,\Delta\big{\vert}_E)$
    is surjective for $\Delta:=C+K_X\big{\vert}_S-2E \in \Pic S$, 
    \label{item:surjective}
    \item $\beta$ is not contained in $H^0(S,N_{S/X})$,
    equivalently, the principal part $\beta\big{\vert}_E$ 
    of $\beta$ is nonzero
    in $H^0(E,N_{S/X}(E)\big{\vert}_E)$,
    \label{item:principal part}
    \item $E$ is an irreducible curve of arithmetic genus of $g(E)$
    and $(\Delta.E)=2g(E)-2-2E^2$, 
    \label{item:numerical condition}
    \item the $\pi$-map $\pi_{E/S}(E)$ is not surjective, and
    \label{item:non-surjective of pi}
    \item $H^1(S,C-E)=0$.
    \label{item:vanishing}
  \end{enumerate}
\end{thm}

Here the natural map in the item \eqref{item:injective} is
induced by an inclusion
$\iota: \mathcal O_S(kE) \hookrightarrow \mathcal O_S((k+1)E)$
of sheaves on $S$.
    
\begin{rmk}
  \label{rmk:obstruction}
  In \cite{Mukai-Nasu,Nasu5} the authors reduced
  the computation of $\ob(\alpha)$ to
  the image $\ob_S(\alpha)$ ($=\pi_{C/S}(\ob(\alpha))$)
  in $H^1(S,N_{S/X}\big{\vert}_C)$.
  They deduced the nonzero of $\ob_S(\alpha)$ 
  from the nonzero of a cup product on the (polar) curve $E$.
  In fact, let $\partial_E$ denote the coboundary map of 
  the short exact sequence
  \begin{equation}
    \label{ses:pi-map}
    \begin{CD}
      [0 @>>> N_{E/S} @>>> N_{E/X} @>{\pi_{E/S}}>> 
	N_{S/X}\big{\vert}_E @>>> 0]\otimes_E \mathcal O_E(E)
    \end{CD}
  \end{equation}
  on $E$. Then the nonzero of $\ob_S(\alpha)$ is deduced from
  that of $\partial_E(\beta\big{\vert}_E)\cup \beta\big\vert_E$,
  where the cup product is taken by the map
  $$
  H^1(E,\mathcal O_E(2E))\times
  H^0(E,N_{S/X}(E-C)\big{\vert}_E) 
  \overset{\cup}\longrightarrow
  H^1(E,N_{S/X}(3E-C)\big{\vert}_E).
  $$
  The condition \eqref{item:numerical condition} assures us that
  the invertible sheaf $N_{S/X}(E-C)\big{\vert}_E$ on $E$ is trivial,
  while \eqref{item:non-surjective of pi} and \eqref{item:vanishing} imply that the coboundary image
  $\partial_E(\beta\big{\vert}_E)$ is nonzero.
  Thus they obtained the nonzero of the cup product on $E$.
  Theorem~\ref{thm:obstrcutedness theorem} looks technical,
  however it has many application in e.g.~\cite{Nasu4,Nasu5,Nasu6}.
  We refer to \cite{Nasu5} for the proof.
\end{rmk}

Now we assume that $X$ is an Enriques-Fano $3$-fold and
$S$ is its hyperplane section. Let $E$ be a half pencil on $S$.
\begin{lem}
  \begin{enumerate}
    \item The $\pi$-map $\pi_{E/S}(E)$ for $(E,S)$ is not surjective
    if and only if $H^1(E,N_{E/X}(E))=0$.
    \item Suppose that there exists a commutative diagram
    \begin{equation}
      \xymatrix{
	F \ar@{^{(}->}[r] 
	\ar[d]_{2:1} 
	& M \ar@{^{(}->}[r] \ar[d]_{2:1}
        & Y \ar[d]^{\pi}_{2:1} \\
	E \ar@{^{(}->}[r] & S \ar@{^{(}->}[r] & X, \\
	  } 
     \end{equation}
    where $\pi$ is the canonical cover of $X$,
    $M$ is the $K3$-cover of $S$ in $Y$
    and $F$ is the pullback of $E$ in $M$.
    Then we have $H^1(E,N_{E/X}(E))=0$ if $H^1(F,N_{F/Y})=0$.
  \end{enumerate}
  \label{lem:not surjective}
\end{lem}

\Proof
(1)\quad 
By adjunction,
we have $H^1(E,N_{S/X}(E)\big{\vert}_E) \simeq 
H^1(E,-K_X\big{\vert}_E +K_E)=0$.
Thus it follows from \eqref{ses:pi-map} that
$\pi_{E/S}(E)$ is surjective if and only if
the map 
$$
\Phi: H^1(E,N_{E/S}(E)) \longrightarrow H^1(E,N_{E/X}(E))
$$
induced by the sheaf homomorphism
$[N_{E/S} \hookrightarrow N_{E/X}]\otimes_E \mathcal O_E(E)$
is an isomorphism.
Since $N_{E/S}(E)\simeq \mathcal O_E(2E)$ is trivial,
we see that $H^1(E,N_{E/S}(E))\simeq k$. 
Thus we obtain the first assertion.

(2)\quad
Note that $X$ is a quotient scheme $Y/G$ of $Y$ 
by a finite group scheme $G$ (of order $2$).
Then $N_{E/X}(E)$ is the $G$-invariant part of $\pi_*N_{F/Y}(F)$.
Therefore $N_{E/X}(E)$ is a direct summand of $\pi_*N_{F/Y}(F)$,
due to the existence of the Reynolds operator
(cf.~\cite[\S2]{Schlessinger71}).
Then there exists a natural injection
$$
H^1(E,N_{E/X}(E)) \hookrightarrow H^1(F,N_{F/Y}(F)).
$$
Since $\mathcal O_F(F)\simeq N_{F/M}$ is trivial,
we conclude that $H^1(E,N_{E/X}(E))=0$ by assumption.
\qed

\section{Deformations of curves on Enriques-Fano $3$-folds}
\label{sect:obstructions}

In this section, we prove Theorem~\ref{thm:main2}.
Let $X$ be an Enriques-Fano $3$-fold of genus $g$,
$S$ an Enriques surface in $X$ and
$C$ a smooth connected curve on $S$ of genus $g(C)$.

\paragraph{\bf Proof of Theorem~\ref{thm:main2}}

We first show a strategy of the proof, which is 
very similar to that of \cite[Theorem~1.2]{Nasu5}.
By Lemma~\ref{lem:nonsingularity of flag},
we have $H^1(X,N_{(C,S)/X})=0$, which implies that
the Hilbert-flag scheme $\HF^{sc} X$ of $X$ is nonsingular at $(C,S)$.
Moreover, it follows from \eqref{ses:normal1} that
there exists an exact sequence
\begin{equation}
  \label{ses:long exact}
  \begin{CD}
    H^0(X,N_{(C,S)/X})
    @>{p_1}>> H^0(C,N_{C/X})
    @>{\delta}>> H^1(S,N_{S/X}(-C))
    @>>> 0,
  \end{CD}
\end{equation}
where $p_1$ is the tangent map of the first projection
$pr_1: \HF^{sc} X \rightarrow \Hilb^{sc} X$, $(C',S') \mapsto [C']$.
We define a divisor $D$ on $S$ as in the statement.
Then since $N_{S/X}(-C)\simeq \mathcal O_S(K_S-D)$,
the cokernel of $p_1$ is isomorphic to $H^{1}(S,D)^{\vee}$ by Serre duality.
Therefore, $p_1$ is surjective if $H^1(S,D)=0$
and we have proved Theorem~\ref{thm:main2}(1) 
by virtue of Lemma~\ref{lem:nonsingularity of flag}.
On the other hand, under the settings of (2) and (3) of the theorem,
by Lemma~\ref{lem:KL.etc}, 
$H^1(S,D)$ is of dimension $m$ and $1$, respectively.
Therefore, there exists a global section $\alpha$ of $N_{C/X}$
not contained in the image of $p_1$.
For the proofs of (2) and (3), it suffices to
prove that the primary obstruction $\ob(\alpha)$ is nonzero
for such an $\alpha$. We prove this only for (2) and 
skip the proof of (3), because in the latter case, 
$S$ is a {\em nodal} Enriques surface
and the proof of $\ob(\alpha)\ne 0$ is more similar to 
that of \cite[Theorem~1.2 (2)]{Nasu5}.
We refer to \cite[\S4]{Nasu5} for more details of the proof of (3).

We prove (2). Suppose that there exist
a half pencil $E$ on $S$ and an integer $m \ge 1$
such that $D\sim 2mE$ or $D \sim K_S+(2m+1)E$.
Then by Lemma~\ref{lem:KL.etc}, we see that
$h^1(S,D)=m$ and $h^1(S,D-E)=m-1$.
Since $E$ is effective, there exists a natural map
$H^1(S,N_{S/X}(-C))\rightarrow H^1(S,N_{S/X}(E-C))$,
where $N_{S/X}(-C)\simeq \mathcal O_S(K_S-D)$ in $\Pic S$.
Then by dimension, there exists a nonzero element $\gamma$ 
in the kernel of this map, i.e.~we have 
$r(\gamma,E)=0$, using the same notation in \eqref{eqn:lifting}.
Let $\delta$ be the coboundary map of \eqref{ses:normal1}.
Then it follows from \eqref{ses:long exact} that
there exists a global section $\alpha$ of $N_{C/X}$
such that $\delta(\alpha)=\gamma$.
Let $\pi_{C/S}(\alpha)$ be the exterior component of $\alpha$
(cf.~\S\ref{subsec:flag and primary obstruction}).
We see that $\delta$ factors through
the coboundary map of the short exact sequence
\begin{equation}
  \label{ses:k_C}
  \begin{CD}
    [0 @>>> \mathcal O_S(-C)
    @>>> \mathcal O_S
    @>{|_C}>> \mathcal O_C
    @>>> 0]\otimes N_{S/X}
\end{CD}
\end{equation}
on $S$ (cf.~\cite[\S2.2]{Nasu6}).
Thus we obtain that $\gamma=\pi_{C/S}(\alpha)\cup \mathbf k_C$
for the extension class
$k_C \in \Ext^1(\mathcal O_C,\mathcal O_S(-C))$ of \eqref{ses:k_C}.
Since the reduction $r(*,E)$ and the cup product map 
$\cup \mathbf k_C$ are compatible, we have
$$
r(\pi_{C/S}(\alpha),E) \cup \mathbf k_C 
=r(\pi_{C/S}(\alpha) \cup \mathbf k_C, E)
=r(\gamma, E)=0.
$$
Then it follows from the exact sequence \eqref{ses:k_C}
tensored with $\mathcal O_S(E)$ that
there exists an element $\beta$ in $H^0(S,N_{S/X}(E))$
(i.e.~an infinitesimal deformation of $S$ with a pole along $E$)
such that $\beta\big{\vert}_C=r(\pi_{C/S}(\alpha),E)$
in $H^0(C,N_{S/X}(E)\big{\vert}_C)$.
It is easy to check that all the conditions (from \eqref{item:injective} to 
\eqref{item:vanishing})
of Theorem~\ref{thm:obstrcutedness theorem} are satisfied.
In fact, \eqref{item:injective} is clear (cf.~\S\ref{subsec:enriques}).
Since $\Delta=C+K_X\big{\vert}_S-2E$,
we have $\Delta \sim (2m-2)E$ or $\Delta \sim K_S+(2m-1)E$
and hence we obtain \eqref{item:surjective}.
Since we have $\pi_{C/S}(\alpha)\cup k_C=\gamma \ne 0$ in $H^1(S,N_{S/X}(-C))$,
\eqref{item:principal part} is a consequence of \cite[Lemma~3.1]{Nasu5}.
\eqref{item:numerical condition} follows from $E^2=0$ and $g(E)=1$.
Since $H^1(E,N_{E/X}(E))=0$,
\eqref{item:non-surjective of pi} follows from
Lemma~\ref{lem:not surjective}.
Since $-K_X$ is ample, so is $C-E=-K_X\big{\vert}_S+E+\Delta$
and hence we obtain \eqref{item:vanishing}.
Thus we have proved (2) of Theorem~\ref{thm:main2}.

Finally we prove the last statement, which is 
concerned with the dimension of $\Hilb^{sc} X$.
Let $\mathcal O_{X,x}$ denote the local ring of a scheme $X$ at a point $x$.
We note that $H^0(S,N_{S/X}(-C))\simeq H^0(S,K_S-D)=0$.
Then it follows from \cite[Theorem~2.4]{Nasu6} that
there exist inequalities
$$
\dim \mathcal O_{\HF^{sc} X,(C,S)}
\le \dim \mathcal O_{\Hilb^{sc} X,[C]}
\le \dim \mathcal O_{\HF^{sc} X,(C,S)}+h^1(S,D),
$$
and the inequality to the right is strict if and only if 
$C$ is obstructed in $X$.
By assumption, we have $h^1(S,D)\le 1$ and
$C$ is obstructed in $X$ if $h^1(S,D)=1$.
Thus we have 
$\dim \mathcal O_{\Hilb^{sc} X,[C]}
=\dim \mathcal O_{\HF^{sc} X,(C,S)}=g+g(C)-1$
by Remark~\ref{rmk:expected dimension}.
\qed

\begin{rmk}
  \label{rmk:not applied}
  Let $D:=C+K_X\big{\vert}_S$.
  If $D \sim (2m+1)E$ or $D \sim K_S+(2m+2)E$ for $m \ge 1$,
  then we still have $H^1(S,D)\ne 0$.
  However, Theorem~\ref{thm:obstrcutedness theorem}
  does not apply to $C$ in this case.
  In fact, in the proof of Theorem~\ref{thm:obstrcutedness theorem},
  the nonzero of $\ob(\alpha)$ is reduced to that of 
  the cup product 
  $\partial_E(\beta\big{\vert}_E)\cup \beta\big{\vert}_E$
  in $H^1(E,N_{S/X}(3E-C)\big{\vert}_E)$
  (cf.~Remark~\ref{rmk:obstruction}).
  We see that this cohomology group is zero in the above case,
  because $N_{S/X}(3E-C)\simeq \mathcal O_S(K_S-D-3E)$
  and we have
  $H^1(E,\mathcal O_E(nE))=0$ for odd $n$ and
  $H^1(E,\mathcal O_E(K_S+nE))=0$ for even $n$.
\end{rmk}

\section{Non-reduced components of the Hilbert scheme}
\label{sect:non-reduced}

In this section, we prove Theorem~\ref{thm:main1}.
We also give some examples of Enriques-Fano $3$-folds 
satisfying the assumption of the theorem
(cf.~Example~\ref{ex:example of GNRC}).
In our examples,
every Enriques-Fano $3$-fold $X$ has only terminal cyclic quotient
singularities and 
there exist a smooth Fano $3$-fold $Y$ that double covers $X$.
We also prove that
$\Hilb^{sc} Y$ also contains a generically non-reduced component
and compare its properties (e.g. the dimension of the component) 
with that of $\Hilb^{sc} X$ (cf.~Remark~\ref{rmk:double dimension}).
In what follows, we fix an Enriques-Fano $3$-fold $X$ of genus $g$
and a smooth hyperplane section $S$ of $X$
and consider (a family of) curves on $X$ (or $S$).

\medskip

\paragraph{\bf Proof of Theorem~\ref{thm:main1}}
We consider a complete linear system
$$
\Lambda:= |-K_X\big{\vert}_S+2E|
$$
of divisors on $S$ for the half pencil $E$ in the theorem.
Let $\Phi$ denote the Cossec-Dolgachev function
(cf.~\cite{Dolgachev16,Cossec-Dolgachev}).
Then we have $\Phi(-K_X\big{\vert}_S+2E)=(-K_X.E)=e\ge 2$,
which indicates that $\Lambda$ is base-point-free
(cf.~\cite[Chap.~IV.~\S4]{Cossec-Dolgachev}).
By Bertini's theorem, $\Lambda$ contains a smooth connected 
member $C$,
which is a curve on $S$ of genus $g(C)=C^2/2+1=g+2e$.
Since $\mathcal O_S(C)$ is ample, we see that $H^1(S,\mathcal O_S(C))=0$.
Then by Lemma~\ref{lem:nonsingularity of flag},
the Hilbert-flag scheme $\HF^{sc} X$ is nonsingular at $(C,S)$
of expected dimension $\chi(X,N_{(C,S)/X})=g+g(C)-1$ 
(cf.~Remark~\ref{rmk:expected dimension}).
Therefore there exists a unique irreducible component $\mathcal W$
of $\HF^{sc} X$ passing through $(C,S)$.
Let $W$ be its image by the first projection 
$pr_1: \HF^{sc} X\rightarrow \Hilb^{sc} X$, $(C,S)\mapsto [C]$.
Then $W$ is an irreducible closed subset of $\Hilb^{sc} X$.
Since $N_{S/X}(-C)\simeq K_S-2E$,
we see that $H^0(S,N_{S/X}(-C))=0$.
This implies that there exists a (Zariski) open neighborhood 
$\mathcal U \subset \mathcal W$ of $(C,S)$, 
the restriction of $pr_1$ to which is an embedding.
Thus we see that $\dim W=\dim \mathcal W=g+g(C)-1$.
Moreover, since $h^1(S,N_{S/X}(-C))=1$,
it follows from \eqref{ses:long exact} that
$$
h^0(C,N_{C/X})=h^0(X,N_{(C,S)/X})+1=g+g(C).
$$
Applying Theorem~\ref{thm:main2} to $C$,
we see that $C$ is obstructed in $X$.
Moreover, since $\dim W$ attains 
the dimension of $\Hilb^{sc} X$ at $[C]$,
$W$ is an irreducible component of $(\Hilb^{sc} X)_{\red}$.

Let $(C',S')$ be a general member of $\mathcal W$.
Then $S'$ is a (smooth) Enriques surface.
Since $H^1(S,\mathcal O_S)=0$, the Picard group of $S$ does not change 
under the smooth deformation of $S$ and hence $\Pic S' \simeq \Pic S$.
Since $H^1(S,\mathcal O_S(E))=0$, the half pencil $E$ is deformed to 
a half pencil $E'$ on $S'$ (of the same degree).
Thereby $C'$ is linearly equivalent to $-K_X\big{\vert}_{S'}+2E'$ for some $E'$. 
By upper semicontinuity, we have $H^1(E',N_{E'/X}\otimes_{E'} N_{E'/S'})=0$.
Again by Theorem~\ref{thm:main2}, $C'$ is obstructed in $X$.
Thus $\Hilb^{sc} X$ is generically singular along $W$.
Since $g+g(C)-1=2g+2e-1$, we have completed the proof.
\qed

\begin{cor}
  \label{cor:EF3 with terminal cyc quot sing}
  Let $X$ be an Enriques-Fano $3$-fold 
  with only terminal cyclic quotient singularities\footnote{
    i.e.~$X$ is one of the $3$-folds listed in Table~\ref{table:EF3}},
  $Y$ the smooth Fano cover of $X$ and
  $\theta$ an involution of $Y$ such that $X \simeq Y/\theta$.
  If there exist
  \begin{enumerate}
    \item a $\theta$-invariant $K3$ surface $M \subset Y$
    not passing through the fixed points of $\theta$, and
    \item an elliptic fibration on $M$ with an $\theta$-anti-invariant
    elliptic parameter $u$ such that $H^1(F,N_{F/Y})=0$ for $F=(u=0)$,
  \end{enumerate}
  and if moreover $f:=(-K_Y.F)\ge 4$,
  then $\Hilb^{sc} X$ contains a generically non-reduced component
  of dimension $2g+f-1$.
\end{cor}
\Proof 
Let $S \subset X$ be the Enriques quotient of $M$ 
(by $\theta_M: =\theta \bigm|_M$),
that is a smooth hyperplane section of $X$.
Then by Lemma~\ref{lem:elliptic fibration},
the image $E$ of $F$ in $S$ is a half-pencil.
We see that $e=(-K_X.E)=f/2\ge 2$
and $H^1(E,N_{E/X}(E))=0$ by Lemma~\ref{lem:not surjective}.
Thus the corollary follows from Theorem~\ref{thm:main1}.
\qed

\begin{ex}
  \label{ex:example of GNRC}
  Suppose that $X$ is one of the Enriques-Fano $3$-folds $X_g$
  of genus $g=13,9,6$ in Example~\ref{ex:EF3}.
  Then by Lemma~\ref{lem:elliptic fibration on K3},
  there exist a $K3$ surface $M \subset Y$
  and an elliptic fibration on $M$ 
  with $\theta$-invariant fiber $F \subset M$
  satisfying the assumption of 
  Corollary~\ref{cor:EF3 with terminal cyc quot sing}.
  Moreover $F$ is a complete intersection in $Y$ 
  of degree $f=(-K_Y.F)=12,8,6$ for $X=X_{13},X_9,X_6$, respectively.
  Therefore $\Hilb^{sc} X$ contains a generically non-reduced component $W$
  with the following properties:
  \begin{enumerate}
    \item every general member $C$ of $W$ is contained
    in an Enriques surface $S \sim_{\mathbb Q} -K_X$, 
    \item $C \sim -K_X\big{\vert}_S+2E$ in $\Pic S$,
    where $E$ is a half pencil on $S$, and 
    \item $h^0(C,N_{C/X})=\dim W+1$
    and $\dim W=37$ if $X=X_{13}$,
    $\dim W=25$ if $X=X_9$ and
    $\dim W=17$ if $X=X_6$.
  \end{enumerate}
\end{ex}

\begin{rmk}
  \label{rmk:generic member}
  In Example~\ref{ex:example of GNRC}
  every general member $C$ of $W$ is contained in 
  a general hyperplane section $S$ of $X \subset \mathbb P^g$
  by Remark ~\ref{rmk:general section}.
\end{rmk}

It may be worthwhile to note that if a smooth Fano $3$-fold $Y$
contains a $K3$ surface and an elliptic curve on the surface,
then under some extra assumptions,
$\Hilb^{sc} Y$ also contains a generically
non-reduced component.

\begin{prop}[{cf.~\cite{Nasu5,Nasu6}}]
  \label{prop:existence of GNRC for smooth Fano}
  Let $Y$ be a smooth Fano $3$-fold anti-canonically embedded
  and $M$ a smooth $K3$ surface in $Y$.
  If there exists an elliptic curve $F$ on $M$ such that
  $$
  H^1(F,N_{F/Y})=H^1(M,\mathcal O_M(F+K_Y\big{\vert}_M))=0,
  $$
  then $\Hilb^{sc} Y$ contains a generically non-reduced component $V$ 
  such that for its general member $C$, we have
  \begin{enumerate}
    \item $C$ is contained in a smooth $K3$ surface $M' \in |-K_Y|$,
    \item $C \sim -K_Y\big{\vert}_{M'}+2F'$ for some elliptic curve
    $F'$ on $M$, and
    \item $\dim V=(-K_Y)^3+2f+2$ and
    $h^0(C,N_{C/Y})=\dim V+1$,
    where $f$ denotes the degree $(-K_Y.F)$ of $F$.
  \end{enumerate}
\end{prop}
\Proof
By adjunction and Serre duality, we see that
$$
H^1(M,N_{M/Y}(-F)))\simeq H^1(M,\mathcal O_M(-K_Y\big{\vert}_M-F))=0.
$$
It follows from \eqref{ses:normal1} that $H^1(Y,N_{(F,M)/Y})=0$.
We consider a complete linear system 
$\Lambda:=|-K_Y\big{\vert}_M+2F|$ on $M$.
Since $\Lambda$ is base-point-free,
there exists a smooth connected curve $C$ on $S$,
whose genus is computed as $g(C)=C^2/2+1=(-K_Y)^3/2+2f+1$.
Since $H^1(M,F)=0$, by virtue of \cite[Lemma~2.12]{Nasu6}\footnote{
There is a typo in this lemma.
The assumption $\car k\ne 0$ is wrong and $\car k=0$ is correct.},
we deduce from $H^1(Y,N_{(F,M)/Y})=0$ that $H^1(Y,N_{(C,M)/Y})=0$.
This implies that
there exists a first order deformation $\tilde M$ of $M$ in $Y$,
to which $C$ does not lift by \cite[Lemma~2.8]{Nasu6}.
Then \cite[Theorem~1.2 and Corollary 1.3]{Nasu5} show that
$C$ is obstructed in $Y$ and moreover,
there exists a generically non-reduced component $V$ 
of $\Hilb^{sc} Y$ passing through $[C]$.
Since $N_{M/Y}(-C)\simeq -2F$, we have $H^0(M,N_{M/Y}(-C))=0$
and then $\dim V$ is equal to
$$
h^0(C,N_{(C,M)/Y})=(-K_Y)^3/2+g(C)+1=(-K_Y)^3+2f+2,
$$
that is the expected dimension of 
the Hilbert-flag scheme $\HF^{sc} Y$ at $(C,M)$.
It follows from $h^1(M,N_{M/Y}(-C))=1$ and \eqref{ses:normal1}
that $h^0(C,N_{C/Y})=h^0(Y,N_{(C,M)/Y})+1$.
Thus the proposition has been proved.
\qed

\begin{rmk}
  \label{rmk:double dimension}
  One can compare Theorem~\ref{thm:main1} with
  Proposition~\ref{prop:existence of GNRC for smooth Fano}.
  It might be also interesting to note that
  in Example~\ref{ex:example of GNRC}
  the $K3$ surface $M$ 
  and the elliptic curve $F$ satisfy
  the assumption of 
  Proposition~\ref{prop:existence of GNRC for smooth Fano}
  (cf.~Remark~\ref{rmk:-K-normal}).
  Therefore the Hilbert scheme $\Hilb^{sc} Y$ of 
  the smooth Fano cover $Y$ of $X$
  contains a generically non-reduced component $V$.
  Moreover, we have
  $$
  \dim V=(-K_Y)^3+2f+2=2((-K_X)^3+2e+1)=2\dim W,
  $$
  where $W$ is the non-reduced component of $\Hilb^{sc} X$
  in Example~\ref{ex:example of GNRC}.
\end{rmk}

For every example $X$ of the Enriques-Fano $3$-folds 
in Example~\ref{ex:example of GNRC}, the smooth Fano cover
$\pi: Y \rightarrow X$ is \'etale in a neighborhood of $C$. 
Thus this example gives us some affirmative evidence to the following question.

\begin{question}
  \label{ques:naive}
  Let $\pi: Y \rightarrow X$ be a finite covering of a projective scheme $X$, $C$ a smooth curve on $X$. 
  Suppose that $\pi$ is \'etale in a neighborhood of $C$. 
  Then is $\pi^{-1}(C)$ obstructed in $Y$ if so is $C$ in $X$?
\end{question} 

We remark that there is no morphism $\Hilb Y \rightarrow \Hilb X$ in general.

\bibliography{mybib}
\bibliographystyle{abbrv}

\end{document}